\documentclass[10pt]{article}

\usepackage{amssymb}
\usepackage{amsfonts}
\usepackage{amsmath}
\usepackage{amsthm}
\usepackage{setspace}
\usepackage{geometry}
\usepackage{shuffle}

\theoremstyle{plain}

\theoremstyle{definition}

\newcommand{\CC}{{\mathbb C}}

\newcommand{\ZZ}{{\mathbb Z}}

\usepackage{hyperref}

\title{Iterated integrals on affine curves}
\author{Martin T. Luu, Albert Schwarz}

\date{}

\begin{document}

\newcommand{\Addresses}{{
\bigskip
\footnotesize

Martin Luu, \textsc{Department of Mathematics, University of California, Davis} \par \nopagebreak \textit{E-mail address:} \texttt{mluu@math.ucdavis.edu}

\medskip

Albert Schwarz, \textsc{Department of Mathematics, University of California, Davis} \par \nopagebreak \textit{E-mail address:} \texttt{schwarz@math.ucdavis.edu}

}}

\maketitle

\begin{abstract}
Motivated by amplitude calculations in string theory we establish basic properties of homotopy invariant iterated integrals on affine curves.
\end{abstract}

$$$$

\hfill \begin{minipage}{0.5\linewidth}
\emph{To Dmitry Fuchs, with warmest wishes for his 80th birthday}
\end{minipage}

\section{Introduction}

In the present paper we consider iterated integrals on a compact two-dimensional smooth manifold with a finite and non-zero number of deleted points. Such a manifold can be considered as an algebraic curve with deleted points or as an affine one-dimensional complex manifold. It is a Stein manifold; this simplifies drastically the description of iterated integrals. The next step will be the analysis of iterated integrals on the configuration spaces on affine curves; these spaces can be considered as affine manifolds.

The configuration spaces on compact complex curves appear naturally in the calculation of string amplitudes.
Iterated integrals can be used in these calculations, however  on compact curves of genus $>0$ it is necessary to work with
non-holomorphic integrals.  

Alternatively  in the calculation of string amplitudes one can work with configuration spaces on affine curves. The results of the present paper should be useful in this approach.

\section{Basic definitions}
Let $X$ be a smooth manifold. Consider smooth $1$-forms $w_{1},\cdots, w_{r}$ on $X$ and a smooth path $\gamma$ between two points $Q_{1}$ to $Q_{2}$ in $X$. Write $\gamma^{*} w_{i} = g_{i} \textrm{ d}t$ and consider the iterated integral
\begin{eqnarray}
\label{LI-definition}
\textrm{L}_{w_{1} \cdots w_{r} }(Q_{1},Q_{2}) &:=& \int_{\gamma} (w_{r} \cdots w_{1}) \\
&& \nonumber\\
&=& \int_{0 \le t_{1} \le \cdots \le t_{r} \le 1} g_{r}(t_{1}) \textrm{ d}t_{1} \cdots g_{1}(t_{r}) \textrm{ d}t_{r} \nonumber
\end{eqnarray}
We will also use the notation 
$$\textrm {L}_w=\int _{\gamma}w^{\textrm {op}}$$ for this integral. Here the superscript $\textrm{op}$ indicates the reversal of the order of the letters in the word $w$. (The integral corresponds to a word $w=w_1\otimes \cdots \otimes w_r$ considered as an element of the tensor power of the space $\Omega^1 (X)$ of $1$-forms). Extend this definition by linearity to all elements of the tensor algebra generated by $\Omega^{1}(X)$. (We will write $w_1\cdots w_r$ instead of $w_1\otimes \cdots\otimes w_r$ and we will use term ``generalized words'' for linear combinations of words). To obtain via Equation (\ref{LI-definition}) a (multi-valued) function of $Q_{2}$ one would like to restrict to those elements where the iterated integral only depends on the homotopy class of $\gamma$. In this case these functions can be seen to generalize the classical polylogarithms. By general results of Chen, see for example \cite{CHE}, the homotopy invariant integrals can be described in terms of the operator $D$ defined by
$$D (w_{1} \otimes \cdots \otimes w_{r})=\sum_{i=1}^{r} w_{1}\otimes \cdots \otimes dw_{i} \otimes \cdots\otimes  w_{r} + \sum_{i=1}^{r-1} w_{1}\otimes \cdots \otimes(w_{i} \wedge w_{i+1})\otimes \cdots  \otimes w_{r}$$
and extended to general elements of the tensor algebra by linearity. Fixing a connected model of the de Rham complex the homotopy functionals correspond to the elements of the tensor algebra generated by  degree $1$ elements of the connected model on which the operator $D$ vanishes. We will establish basic properties of the corresponding multi-valued functions obtained via Equation (\ref{LI-definition}). We write from now on $z$ for $Q_{2}$ and we suppress the choice of $Q_{1}$ in the notation and denote the multi-valued functions of Equation (\ref{LI-definition}) by $\textrm{L}_{w}(z)$. If $X$ is an affine manifold then we can restrict ourselves to $w$'s that are generalized words with respect to holomorphic forms; see \cite{HAI} (Section 13) for a detailed discussion. 

In the present paper we consider the case where $X$ is obtained from a compact smooth surface $X'$ by deleting a finite subset $S=\{P_{1},\cdots,P_{n}\}$. We assume that the set of deleted points is not empty; then $X$ can be realized as an affine algebraic curve and the homotopy invariant iterated integrals can be very simply described. They correspond to words with respect to a basis of one-dimensional cohomology represented by holomorphic $1$-forms on $X$. (Such a word represents a homotopy invariant iterated integral because on a curve every holomorphic $1$-form is closed and a wedge product of two holomorphic $1$-forms vanishes. From the other side for affine curves one can find a basis in one-dimensional cohomology represented by holomorphic $1$-forms).

The dimension of the first de Rham cohomology of $X=X'\backslash \{P_{1},\cdots,P_{n}\}$ is known to equal $2g+n-1$ where $g$ denotes the genus of the surface $X'$. One can choose a basis such that $g+n-1$ of the basis elements are logarithmic $1$-forms. In particular, on the punctured sphere the logarithmic $1$-forms constitute the whole basis. Namely, if $z$ denotes a coordinate on $X'\backslash \{P_{n}\}$ then one can take as a basis the forms $w_{i}=(z-P_{i+1})^{-1} \textrm{ d}z$ for $0\le i \le n-2$. In genus $0$ foundational results on iterated integrals and their relation to multiple zeta values are described in \cite{BRO} and \cite{GON}. We show that many results carry over to iterated integrals on surfaces of arbitrary genus. 

We consider in some detail iterated integrals on surfaces of genus 1 (on affine elliptic curves) in the holomorphic picture (iterated integrals on elliptic curves were considered by Brown and Levin in \cite {BL} in the non-holomorphic picture). Let us suppose that the elliptic curve $X'$ is represented as $\mathbb {C}/(\ZZ+\tau \ZZ)$ for a complex number $\tau$ with positive imaginary part. A basis of the first de Rham cohomology of $X'\backslash \{P_{1},\cdots,P_{n}\}$ can be chosen to consist of logarithmic $1$-forms $w_0,\cdots , w_{n-1}$ and a $1$-form $w_n$ having a second order pole at $P_{1}$, say. One can assume that $w_0$ is holomorphic and that $w_i$ with $1\le  i\le  n-1$ has simple poles at $P_{i+1}, P_{1}$. 

In Section \ref{zeta-section} an important role is played by the shuffle product relations that the functions $\textrm{L}_{w}(z)$ satisfy. Let $\shuffle$ denote the shuffle product in the tensor algebra. It is defined as
$$w_{a_{1}}\cdots w_{a_{r}} \shuffle w_{a_{r+1}} \cdots w_{a_{s}} := \sum_{\sigma} w_{a_{\sigma(1)}} \cdots w_{a_{\sigma(r+s)}}$$
where $\sigma$ ranges over those permutations on $r+s$ elements such that $\sigma^{-1}(i) < \sigma^{-1}(j)$ whenever $i<j$ and either $i$ and $j$ are both at most $r$ or both at least $r+1$. (Again we extend this definition to arbitrary elements of the tensor algebra by linearity). It is easy to check the shuffle product relation
\begin{eqnarray}
\label{shuffle-relation}
\textrm{L}_{w_{a_{1}}\cdots w_{a_{r}}}(z)\cdot \textrm{L}_{w_{a_{r+1}}\cdots w_{a_{s}}}(z) = \textrm{L}_{w_{a_{1}}\cdots w_{a_{r}} \shuffle w_{a_{r+1}} \cdots w_{a_{s}}}(z)
\end{eqnarray}

\section{Multiple zeta values}
\label{zeta-section}

We consider special values of the functions $\textrm{L}_{w}(z)$ where as before we let $X$ be an affine curve. Due to the relation with classical polylogarithms we denote in this  situation the functions by $\textrm{Li}_{w}(z)$. The special values that we consider generalize the classical (multiple) zeta values and are defined in the following manner: Instead of choosing the starting point $Q_{1}$ and endpoint $Q_{2}$ of the smooth path $\gamma$ to lie in $X=X'\backslash S$ we assume $Q_{1},Q_{2}$ are in the set $S$ of deleted points. Of course, in order to obtain well defined expressions one has to choose in general a suitable regularization procedure for the iterated integrals between two punctures.

The relation between zeta functions and homotopy invariant iterated integrals between punctures on a surface can be seen from the following special case: Consider the sphere $X'=\mathbb{P}^{1}$ and $S=\{0,1,\infty\}$ and let $Q_{1}=0$ and $Q_{2}=1$. Let $z$ denote a coordinate on $X'\backslash \{\infty\}$ and take the basis of $1$-forms on $X=X'\backslash S$ as $w_{0} = \textrm{d}z/z$ and $w_{1} = \textrm{d}z/(1-z)$. Then the direct calculation of the iterated integral gives for $n\ge 2$
$$\textrm{Li}_{w_{0}^{n} w_{1}}(Q_{1},Q_{2}) = \sum_{k=1}^{\infty}\frac{1}{k^{n}}=\zeta(n)$$
Hence special values of the zeta function are expressed as special iterated integrals where the path $\gamma$ connects two punctures. 

\subsection{Regularization near a puncture}

We now describe the shuffle regularization of iterated integrals in the case where the starting point of the path $\gamma$ is taken to be a puncture $P_{j}$. We fix a basis of the first cohomology consisting of holomorphic forms on the affine curve (on the punctured surface); we assume that only one of these forms (call it $w_{j}$) has a pole at $P_{j}$ and it is simple. (It seems that the uniqueness of the form $w_{j}$ is not important, but the simplicity of the pole is). We call the collection of punctures with this property the set of good punctures. For example, in genus $1$ the basis of $1$-forms can be chosen so that all punctures except one are good.

In the case of a good puncture we now show that one can apply the known genus $0$ regularization arguments. Let as before $z$ be in $X=X'\backslash \{P_{1},\cdots, P_{n}\}$. Fix a smooth path from $P_{j}$ to $z$. For each $k$ there is a unique constant $a_{k}$ (zero unless $k=j$) such that 
$$\lim_{\epsilon \rightarrow 0} \; \; \;  \left (-a_{k} \log \epsilon +\int_{\gamma_{\epsilon}} w_{k} \right )$$
exists, where $\gamma_{\epsilon}$ is the piece of $\gamma$ 
between $\epsilon$ and $z$ (we assume here that the coordinate of the puncture is equal to zero and $z$ denotes both the point $Q_2$ and its coordinate). Define $\int_{\gamma} w_{k}$ to be the value of this limit. More generally, for each positive integer $r$ define
$\int_{\gamma} w_{k}^{r}= ( \int_{\gamma} w_{k} )^{r}/r!$ (this definition is prompted by the shuffle product relation).
 
We use the shuffle product to define more generally the value of the iterated integral $\int_{\gamma} w$ where $w$ is an arbitrary generalized word with respect to $w_{0},\cdots,w_{n-1}$. We claim that every word $w$ with respect to $w_{0},\cdots, w_{n-1}$ has a unique expression of the form
\begin{eqnarray}
\label{decomposition-equation}
w= \sum_{i=0}^{k}  w(i) \shuffle w_{j}^{i}
\end{eqnarray}
where $w(i)$ is a word not ending in $w_{j}$. To show existence, consider the summand $vw_{j}^{k}$ of $w$ ending in a maximal amount of $k$ copies of $w_{j}$'s. Subtracting $v \shuffle w_{j}^{k} $ from $w$ yields a word where every summand ends in at most $k-1$ $w_{j}$'s and the desired existence follows by induction on $k$. To show uniqueness suppose there are two decompositions as in Equation (\ref{decomposition-equation}) 
$$ \sum_{i=0}^{k}   w(i) \shuffle w_{j}^{i}=w= \sum_{i=0}^{l} w'(i) \shuffle w_{j}^{i} $$
Then the summand of $w$ ending in the largest amount of $w_{j}$'s is given by $w(k)w_{j}^{k}$ and $w'(l)w_{j}^{l}$ and hence $k=l$ and $w(k)=w'(k)$ and the uniqueness follows by induction on $k$. 

Suppose now a decomposition as in Equation (\ref{decomposition-equation}) is given. We have already defined
$$\textrm{Li}_{w_{j}^{i},j}(z):= \int_{\gamma} w_{j}^{i}=\left (\int_{\gamma} w_{j} \right )^{i}/i! $$
Since $P_{j}$ is a good puncture and since $w(i)$ does not end with $w_{j}$ one can define
$$\textrm{Li}_{w(i),j}(z) :=  \int_{\gamma} w(i)^{\textrm{op}}$$
Prompted by the shuffle product relations we define
\begin{eqnarray}
\label{complete-regularization-equation}
\textrm{Li}_{w,j}(z) := \sum_{i=0}^{k} \textrm{Li}_{w_{j}^{i},j}(z)\textrm{Li}_{w(i),j}(z) 
\end{eqnarray}
This is our choice of regularization of the functions from Equation (\ref{LI-definition}) where the starting point of the path $\gamma$ is a good puncture.

Notice that in our definition of regularization we have imposed some conditions on the puncture. Instead of these conditions we could impose the requirement that the words under consideration consist of logarithmic forms.

\subsection{Multiple zeta values and monodromy}

The results of the previous section can be used to define multiple zeta 
values as iterated integrals along paths between two deleted points. For 
the remainder of this section fix two good punctures $P_{j}$ and $P_{i}$ and a smooth path $\gamma$ from $P_{j}$ to $P_{i}$. Under this assumption the genus $0$ calculations given in \cite{BRO} can be generalized to the current situation. Below we describe this in some detail.

For each generalized word $w$ there is a non-negative integer $t$ such that the function $\textrm{Li}_{w,j}(z)$ defined in Equation (\ref{complete-regularization-equation}) has the asymptotic behavior
\begin{eqnarray}
\label{aymptotic-equation}
\textrm{Li}_{w,j}(z) \sim \sum_{s=0}^{t} a_{s}(z) \log(P_{i}-z)^{s}
\end{eqnarray}
as $z$ approaches $P_{i}$. Here the functions $a_{s}(z)$ have a well defined limit as $z \rightarrow P_{i}$. Define the  multiple zeta value as
$$\textrm{MZV}_{i,j}(w):= \lim_{z \rightarrow P_{i}} \;\; a_{0}(z)$$

We now show that the multiple zeta values can also be realized as the coefficients of a suitable analogue of Drinfeld associators. Let ${\bf x}_{0}, \textbf{x}_{1}, \cdots $ be free non-commuting variables. For a good puncture $P_{k}$ define the generating function 
$$\textrm{L}_{k}(z):= \sum_{w} \textrm{Li}_{w,k}(z) \; {\bf x} $$
where the summation is over all words $w$ with respect to $w_{0}, w_{2},\cdots, w_{r}$ and $\textbf{x}$ is the corresponding word with respect to ${\bf x}_{0}, \textbf{x}_{1}, \cdots $. By convention the coefficient of the empty word is $1$. This generating function is a formal power series in free noncommuting variables ${\bf x}_{0}, \textbf{x}_{1}, \cdots $ whose coefficients are functions of $z$. We now show that $\textrm{L}_{j}(z)$ as well as $\textrm{L}_{i}(z)$ satisfy the same differential equation; this allows us to relate their quotient to multiple zeta values.

First note that for two paths $\alpha$ and $\beta$ in $X\backslash S$ with the starting point of $\beta$ equal to the end point of $\alpha$,
we can say that
$$\int_{\beta \circ \alpha} w_{r}  \cdots w_{1} - \int_{\alpha} w_{r}  \cdots w_{1}  = \int_{\alpha}w_{r} \cdots w_{2} \int_{\beta} w_{1} + \int_{\alpha} w_{r} \cdots w_{3}  \int_{\beta} w_{2} w_{1}
+ \cdots$$
(See for example \cite{BRO} (Prop. 2.2 (iii))).

To calculate $d \; \textrm{Li}_{w}(z)$ let $\alpha$ be 
the path from $Q_{1}$ to $z$ and $\beta$ a path from $z$ to $z+h$, 
for $h$ small. Using that there is a constant $C$ such that the iterated 
integral satisfies
$$|\int_{\gamma} w_{s} \cdots w_{1} | \le C\cdot \textrm{length}
(\gamma)^{s}$$
one obtains the formula
\begin{eqnarray}
d \; \textrm{Li}_{w_{i_{1}} \cdots w_{i_{s}}}(z) = w_{i_{1}} \; \textrm{Li}_{w_{i_{2}}\cdots w_{i_{s}}}(z)
\end{eqnarray}
The above iterated integrals are integrals over a path between two points in $X'\backslash \{P_{1},\cdots,P_{n}\}$ but from our regularization scheme one can see that the same equation holds if the starting point of the integration path is a good puncture $P_{k}$:
\begin{eqnarray}
\label{z-differentiation}
d \; \textrm{Li}_{w_{i_{1}} \cdots w_{i_{s}},k}(z) = w_{i_{1}} \; \textrm{Li}_{w_{i_{2}}\cdots w_{i_{s}},k}(z) 
\end{eqnarray}
From Equation (\ref{z-differentiation}) it follows that
\begin{eqnarray*}
d \; \textrm{L}_{k}(z)= \sum_{w} d\; \textrm{Li}_{w_{i_{1}}\cdots w_{i_{s}},k}(z) \; \textbf{x}_{i_{1}} \cdots \textbf{x}_{i_{s}} &=& \sum_{w} w_{i_{1}} \textbf{x}_{i_{1}} \; \textrm{Li}_{w_{i_{2}}\cdots w_{i_{s}},k}(z) \; \textbf{x}_{i_{2}} \cdots \textbf{x}_{i_{s}}\\
&&\\
&=&\sum_{t=0}^{r} w_{t} \textbf{x}_{t} \cdot \textrm{L}_{k}(z)
\end{eqnarray*}
hence
\begin{eqnarray}
\label{same-differential-equation}
\left (d - \sum_{t=0}^{r} w_{t} {\bf x}_{t} \right )\; \textrm{L}_{j}(z) = 0= \left (d - \sum_{t=0}^{r} w_{t} {\bf x}_{t} \right )\; \textrm{L}_{i}(z)
\end{eqnarray}
The generating series $\textrm{L}_{i}(z)$ has an inverse in the space of formal power series with respect to $\textbf{x}_{0},\textbf{x}_{1},\cdots$ (the coefficients are multivalued functions of  $z$) and we define 
\begin{eqnarray}
\label{Phi-definition}
\Phi_{i,j} = \textrm{L}_{i}(z)^{-1} \cdot \textrm{L}_{j}(z) = \sum_{w} \Phi_{i,j}(w)\; {\bf x} 
\end{eqnarray}
where again $\textbf{x}$ corresponds to a word $w$. Taking the derivative with respect to $z$ and using Equation (\ref{same-differential-equation}) we obtain that the coefficients $\Phi_{i,j}(w)$ are independent of $z$. Therefore $\Phi_{i,j}(w)$ can be calculated by letting $z$ approach the puncture $P_{i}$. As $z$ approaches $P_{i}$ only the iterated integrals of $w_{i}^{k}$ contribute to $\textrm{L}_{i}(z)$. After scaling we can assume that the residue of $w_{i}$ at $P_{i}$ is equal to $1$. Then one has the asymptotic behavior 
\begin{eqnarray}
\label{ijequation}
\textrm{L}_{i}(z)^{-1} \sim \left (\sum_{k=0}^{\infty} \textrm{Li}_{w_{i}^{k},i}(z) \; \textbf{x}_{i}^{k}\right )^{-1} = \left (\sum_{k=0}^{\infty} \frac{\textrm{Li}_{w_{i},i}(z)^{k}}{k!} \; \textbf{x}_{i}^{k}  \right )^{-1} = \exp \left ( - \textrm{Li}_{w_{i},i}(z) \; \textbf{x}_{i}\right )
\end{eqnarray}
and hence 
\begin{eqnarray}
\label{asymptotic-equation}
\textrm{L}_{i}(z)^{-1} \sim \exp(- \log(P_{i}-z) \; \textbf{x}_{i})=\sum_{k=0}^{\infty} (-1)^{k} \cdot  \frac{\log^{k} (P_{i}-z)}{k!} \cdot\textbf{x}_{i}^{k} 
\end{eqnarray} 
One can say also that in the neighborhood of a puncture only one term is relevant in the sum $\sum_{t=0}^{r} w_{t} \textbf{x}_{t}$ in Equation (\ref{same-differential-equation}); neglecting all other terms and solving this equation we obtain the same asymptotic behavior.

Consider now an arbitrary word $w$ and write it as $w = w_{i}^{a} \; w'$ with $a\ge 0$ and $w'$ not beginning with $w_{i}$. Replacing $\textrm{L}_{i}(z)^{-1}$ by its asymptotic expression via Equation (\ref{asymptotic-equation}) one sees that the coefficient $\Phi_{i,j}(w)$ has the following  asymptotic behavior as $z$ goes to $P_{i}$
$$\Phi_{i,j}(w)= \Phi_{i,j}(w_{i}^{a} \; w') \sim \sum_{k=0}^{a} ( \textrm{coefficient of $\textbf{x}_{i}^{k}$ in $ \exp \left ( - \log(P_{i}-z)  \; \textbf{x}_{i} \right )$} )\cdot  ( \textrm{coefficient of $\textbf{x}_{i}^{a-k} \textbf{x}'$ in $\textrm{L}_{j}(z)$} )  $$ 
where $\textbf{x}'$ is the word with respect to $\textbf{x}_{0}, \textbf{x}_{1}, \cdots$ corresponding to $w'$.
It follows that
\begin{eqnarray}
\Phi_{i,j}(w) \sim \sum_{k=0}^{a}  (-1)^{k} \cdot  \frac{\log^{k} (P_{i}-z)}{k!} \cdot  \textrm{Li}_{w_{i}^{a-k} w' ,  j}(z) 
\end{eqnarray}
Since $\Phi_{i,j}$ does not depend on $z$ one has
\begin{eqnarray}
\label{associator-equation}
\Phi_{i,j}(w) = \lim_{z \rightarrow P_{i}} \;\;\;
\sum_{k=0}^{a} (-1)^{k} \cdot  \frac{\log^{k} (P_{i}-z)}{k!} \cdot \textrm{Li}_{w_{i}^{a-k}w',  j}(z) 
\end{eqnarray}
Each function $\textrm{Li}_{w_{i}^{a-k}w' ,  j}(z)$ has an expansion of the form
$$\textrm{Li}_{w_{i}^{a-k}w' ,  j}(z)  = c_{k,0} + c_{k,1}\log(P_{i}-z) + c_{k,2} \log(P_{i}-z)^{2}+ \cdots$$
such that the coefficients $c_{k,0},\cdots$ have a well defined limit as $z \rightarrow P_{i}$. It follows that the only non-singular term of the expression on the right-hand side of Equation (\ref{associator-equation}) is $c_{0,0}$. But this is precisely the term $a_{0}(z)$ associated to $\textrm{Li}_{w_{i}^{a}w',j}(z)$ as in Equation (\ref{aymptotic-equation}). It follows that
$$\Phi_{i,j}(w) = \textrm{MZV}_{i,j}(w)$$
and therefore the coefficients of $\Phi_{i,j}$ are multiple zeta values.

Note also that one can relate the monodromy of the functions $\textrm{Li}_{w,j}(z)$ to the multiple zeta values. Namely, let $\mathcal M_{i}$ be the analytic continuation operator for a loop around $P_{i}$. Since the coefficients of the $\Phi_{i,j}$ are constant it follows that 
$$\mathcal M_{i}(\textrm{L}_{j}(z))= \mathcal M_{i}(\textrm{L}_{i}(z)) \Phi_{i,j}$$
Since we have shown that the coefficients of the $\Phi_{i,j}$ are multiple zeta values this gives an expression of the monodromy around $P_{i}$ in terms of multiple zeta values.

\section{Variational results}
It is important to understand the dependence of the functions $\textrm{Li}_{w}(z)$ on the set of deleted points $S=\{P_{1},\cdots,P_{n}\}$. As indicated before, if the compact surface $X'$ is of genus $0$ and $z$ denotes a coordinate on $X'\backslash \{P_{n}\}$ then one can take as a basis of the first de Rham cohomology the $1$-forms $w_{i}=\mathfrak f_{i}\textrm{ d}z=(z-P_{i})^{-1} \textrm{ d}z$ for $1 \le i \le n-1$. Using the relation
\begin{eqnarray}
\label{partial-fractions-equation}
\mathfrak f_{a}\cdot \mathfrak f_{b} = \sum_{i=1}^{n-1} C_{a,b}^{(i)} \cdot  \mathfrak f_{i}
\end{eqnarray}
where $C_{a,b}^{(i)}=0$ unless $i=a$ or $i=b$ and
 $$C_{a,b}^{(a)} = \frac{1}{P_{a}-P_{b}} \;\;\; ,\;\;\;  C_{a,b}^{(b)}= \frac{1}{P_{b}-P_{a}}$$
one can calculate the action of the operators $\partial_{P_i}$ on the functions $\textrm{Li}_{w}(z)$. If $w_{j_{1}},\cdots, w_{j_{r}}$ is a collection of distinct $1$-forms and $k\ne 1,r$ (the formulas are similar but slightly different for $k=1$ and $k=r$) then  
\begin{eqnarray}
\label{sphere-differential-equation}
& \partial_{P_{j_{k}}}  \; \textrm{Li}_{w_{j_{1}}\cdots w_{j_{k}} \cdots w_{j_{r}}}(z)  \nonumber \\
& =\\
&  \sum_{i=1}^{n-1}  \left ( C_{j_{k},j_{k+1}}^{(i)} \cdot  \textrm{Li}_{w_{j_{1}}\cdots \widehat{w}_{j_{k}} \widehat{w}_{j_{k+1}}  w_{i} \cdots w_{j_{r}}}(z)  -C_{j_{k-1},j_{k}}^{(i)} \cdot  \textrm{Li}_{w_{j_{1}}\cdots \widehat{w}_{j_{k-1}} \widehat{w}_{j_{k}} w_{i} \cdots w_{j_{r}}}(z) \right ) \nonumber
\end{eqnarray}
where the superscript $\; \widehat{} \; $ denotes that the indicated $1$-form is omitted. See for example \cite{GON} for more details.
Since only $C_{a,b}^{(a)}$ and $C_{a,b}^{(b)}$ are non-zero this simplifies to
\begin{eqnarray*}
\partial_{P_{j_{k}}} \; \textrm{Li}_{w_{j_{1} \cdots w_{j_{k}} \cdots w_{j_{r}}}}(z) &=&   C_{j_{k},j_{k+1}}^{(k)} \textrm{Li}_{w_{j_{1}}\cdots\widehat{w}_{j_{k+1}} \cdots w_{j_{r}}}(z) + C_{j_{k},j_{k+1}}^{(k+1)} \textrm{Li}_{w_{j_{1}}\cdots\widehat{w}_{j_{k}} \cdots w_{j_{r}}}(z) \\
&&\\
&&- C_{j_{k-1},j_{k}}^{(k-1)} \textrm{Li}_{w_{j_{1}}\cdots\widehat{w}_{j_{k}} \cdots w_{j_{r}}}(z)- C_{j_{k-1},j_{k}}^{(k)} \textrm{Li}_{w_{j_{1}}\cdots\widehat{w}_{j_{k-1}} \cdots w_{j_{r}}}(z) \\
&&\\
&=& \frac{1}{P_{j_{k-1}} -P_{j_{k}}} \cdot \textrm{Li}_{w_{j_{1} \cdots \widehat{w}_{j_{k-1}} \cdots w_{j_{r}}}}(z)+\frac{1}{P_{j_{k}} -P_{j_{k+1}}} \cdot \textrm{Li}_{w_{j_{1} \cdots \widehat{w}_{j_{k+1}} \cdots w_{j_{r}}}}(z) \\
&& \\
&& -\frac{P_{j_{k-1}}-P_{j_{k+1}}}{(P_{j_{k-1}} -P_{j_{k}})(P_{j_{k}} -P_{j_{k+1}})} \cdot \textrm{Li}_{w_{j_{1} \cdots \widehat{w}_{j_{k}} \cdots w_{j_{r}}}}(z)
\end{eqnarray*}
We now give a generalization to the case where $X'$ is a torus, see \cite{BMMS} for related calculations. To describe the logarithmic $1$-forms fix a complex structure on $X'$, say $X'=\CC/(\ZZ+\tau \ZZ)$ for a complex number $\tau$ with positive imaginary part. Let $f(z)$ be an elliptic function with simple poles precisely at the lattice points $\ZZ+\tau\ZZ$. One can describe $f$ in terms of theta functions, for example one can set
$$f(z)= \partial_{z} \log \theta_{1,1}(z) \;\;\; , \;\;\; \theta_{1,1}(z) = \sum_{n\in \ZZ} \exp \left (\pi i (n+\frac{1}{2})^{2} \tau +2 \pi i (n+\frac{1}{2})(z+\frac{1}{2}) \right )$$  
Then for the space of logarithmic $1$-forms on $X'\backslash \{P_{1},\cdots, P_{n}\}$ one can take the basis $\{w_{0},\cdots,w_{n-1}\}$ where
$$w_{0}= \textrm{d}z = \mathfrak f_{0} \textrm{ d}z $$ 
and
\begin{eqnarray}
\label{elliptic-differentials}
w_{k}=\mathfrak f_{k} \textrm{ d}z  = \left ( f(z-P_{k_{1}})-f(z-P_{k_{2}}) \right ) \textrm{ d}z \;\;\; (1\le k \le n-1)
\end{eqnarray}
where the points $P_{k_{1}},P_{k_{2}}$ in $S$ are suitably chosen.  For example one can take $P_{k_{2}} = P_{1}$ for all $k$ and $P_{k_{1}}=P_{k+1}$. Another possibility is to choose $P_{k_{2}}=P_{k+1}$ and $P_{k_{1}}=P_{k}$ for all $k$. We now consider the iterated integrals in Equation (\ref{LI-definition}) where the $w_{j_{i}}$'s are chosen such that no two of them have poles in common. It follows that for $a\ne b$ there are unique constants $C_{a,b}^{(i)}$ (meaning independent of $z$) such that the elliptic functions $\mathfrak f_{i}$ of Equation (\ref{elliptic-differentials}) satisfy
\begin{eqnarray}
\label{structure-constants-equation}
\mathfrak f_{a}\cdot \mathfrak f_{b} = \sum_{i=0}^{n-1} C_{a,b}^{(i)} \cdot  \mathfrak f_{i}
\end{eqnarray}
Such a decomposition exists since the left-hand side is an elliptic function with simple poles. The elliptic structure constants $C_{a,b}^{(i)}$ can be calculated using the Fay identities for theta functions, see the work of Broedel, Mafra, Matthes, Schlotterer \cite{BMMS} for related calculations. For an integer $\nu$ between $1$ and $n-1$ write $z_{\nu}=z-P_{\nu}$. One can use the Fay identities to show that for distinct integers $i$ and $j$ between $1$ and $n-1$ one has
\begin{eqnarray}
\label{Fay-equation}
\partial_{z} \log \theta_{1,1}(z_{i}) \cdot \partial_{z} \log \theta_{1,1}(z_{j}) &=& \partial_{z} \log \theta_{1,1}(z_{i})\cdot \partial_{z} \log \theta_{1,1}(P_{i}-P_{j}) +\partial_{z} \log \theta_{1,1}(z_{j}) \cdot \partial_{z} \log \theta_{1,1}(P_{j}-P_{i})\\
&& \nonumber \\
&& +\frac{(\partial_{z} \log \theta_{1,1}(z_{j}))^{2}+\partial_{z}^{2} \log \theta_{1,1}(z_{j})}{2}+\frac{(\partial_{z} \log \theta_{1,1}(z_{i}))^{2}+\partial_{z}^{2} \log \theta_{1,1}(z_{i})}{2} \nonumber \\
&& \nonumber \\
&&+\frac{(\partial_{z} \log \theta_{1,1}(P_{i}-P_{j})))^{2}+\partial_{z}^{2} \log \theta_{1,1}(P_{i}-P_{j})}{2} -\frac{1}{2}\cdot \frac{\theta_{1,1}'''(0)}{\theta_{1,1}'(0)} \nonumber
\end{eqnarray}
Fix distinct integers $a$ and $b$ between $1$ and $n-1$. Using Equation (\ref{Fay-equation}) one can show that
$$\mathfrak f_{a}\cdot \mathfrak f_{b}= (\partial_{z} \log \theta_{1,1}(z_{a_{1}}) -\partial_{z} \log \theta_{1,1}(z_{a_{2}})) \cdot  (\partial_{z} \log \theta_{1,1}(z_{b_{1}}) -\partial_{z} \log \theta_{1,1}(z_{b_{2}})) $$
can be written as
\begin{eqnarray*}
 &&\partial_{z} \log \theta_{1,1}(z_{a_{1}})  (\partial_{z} \log \theta_{1,1} (P_{a_{1}}-P_{b_{1}})-\partial_{z} \log \theta_{1,1}(P_{a_{1}}-P_{b_{2}})) \\
&&\\
&& + \partial_{z} \log \theta_{1,1}(z_{a_{2}})  (-\partial_{z} \log \theta_{1,1}(P_{a_{2}}-P_{b_{1}})+\partial_{z} \log \theta_{1,1}(P_{a_{2}}-P_{b_{2}})) \\
&&\\
&&+ \;  \partial_{z} \log \theta_{1,1}(z_{b_{1}})  (\partial_{z} \log \theta_{1,1}(P_{b_{1}}-P_{a_{1}})-\partial_{z} \log \theta_{1,1}(P_{b_{1}}-P_{a_{2}}))\\
&&\\
&&+ \partial_{z} \log \theta_{1,1}(z_{b_{2}})  (-\partial_{z} \log \theta_{1,1}(P_{b_{2}}-P_{a_{1}})+\partial_{z} \log \theta_{1,1}(P_{b_{2}}-P_{a_{2}})) \\
&&\\
&&\;  + \frac{(\partial_{z} \log \theta_{1,1}(P_{a_{1}}-P_{b_{1}})))^{2}+\partial_{z}^{2} \log \theta_{1,1}(P_{a_{1}}-P_{b_{1}})}{2}- \frac{(\partial_{z} \log \theta_{1,1}(P_{a_{1}}-P_{b_{2}})))^{2}+\partial_{z}^{2} \log \theta_{1,1}(P_{a_{1}}-P_{b_{2}})}{2}\\
&&\\
&& \; - \frac{(\partial_{z} \log \theta_{1,1}(P_{a_{2}}-P_{b_{1}})))^{2}+\partial_{z}^{2} \log \theta_{1,1}(P_{a_{2}}-P_{b_{1}})}{2}+\frac{(\partial_{z} \log \theta_{1,1}(P_{a_{2}}-P_{b_{2}})))^{2}+\partial_{z}^{2} \log \theta_{1,1}(P_{a_{2}}-P_{b_{2}})}{2}
\end{eqnarray*}
We suppose now that the pole structure of the basis $w_{0},\cdots, w_{n-1}$ of logarithmic $1$-forms is chosen so that there is an index $i(a,b)$ such that $w_{i(a,b)} = \mathfrak f_{i(a,b)} \textrm{ d}z$ with
$$\mathfrak f_{i(a,b)} =\partial_{z} \log \theta_{1,1}(z-P_{a_{2}}) - \partial_{z} \log \theta_{1,1}(z-P_{b_{2}})$$
One can use this to calculate the elliptic structure constants. One obtains
\begin{eqnarray}
\label{product-equation}
\mathfrak f_{a}\cdot \mathfrak f_{b} = C^{(0)}_{a,b} \cdot  \mathfrak f_{0} + C^{(a)}_{a,b} \cdot\mathfrak f_{a}+C^{(b)}_{a,b}  \cdot \mathfrak f_{b}+C^{(i(a,b))}_{a,b}  \cdot \mathfrak f_{i(a,b)} 
\end{eqnarray}
where the constant coefficients are given by 

\begin{eqnarray*}
C^{(0)}_{a,b} &=&  \frac{(\partial_{z} \log \theta_{1,1}(P_{a_{1}}-P_{b_{1}})))^{2}+\partial_{z}^{2} \log \theta_{1,1}(P_{a_{1}}-P_{b_{1}})}{2}- \frac{(\partial_{z} \log \theta_{1,1}(P_{a_{1}}-P_{b_{2}})))^{2}+\partial_{z}^{2} \log \theta_{1,1}(P_{a_{1}}-P_{b_{2}})}{2}\\
&&\\
&& \; - \frac{(\partial_{z} \log \theta_{1,1}(P_{a_{2}}-P_{b_{1}})))^{2}+\partial_{z}^{2} \log \theta_{1,1}(P_{a_{2}}-P_{b_{1}})}{2}+\frac{(\partial_{z} \log \theta_{1,1}(P_{a_{2}}-P_{b_{2}})))^{2}+\partial_{z}^{2} \log \theta_{1,1}(P_{a_{2}}-P_{b_{2}})}{2}
\end{eqnarray*}
$$C^{(a)}_{a,b} = \partial_{z} \log \theta_{1,1}(P_{a_{1}}-P_{b_{1}})-\partial_{z} \log \theta_{1,1}(P_{a_{1}}-P_{b_{2}})$$
$$C^{(b)}_{a,b} = \partial_{z} \log \theta_{1,1}(P_{b_{1}}-P_{a_{1}})-\partial_{z} \log \theta_{1,1} (P_{b_{1}}-P_{a_{2}})$$
$$C^{(i(a,b))}_{a,b} =\partial_{z} \log \theta_{1,1}(P_{a_{1}}-P_{b_{1}})-\partial_{z} \log \theta_{1,1}(P_{a_{1}}-P_{b_{2}})-\partial_{z} \log \theta_{1,1}(P_{a_{2}}-P_{b_{1}})+\partial_{z} \log \theta_{1,1}(P_{a_{2}}-P_{b_{2}})
$$
In particular, one sees that as in the case of iterated integrals on the sphere, most structure constants vanish. 

We now calculate the variation of the functions $\textrm{Li}_{w}(z)$ as the deleted points move. We assume $k\ne 1,r$ for simplicity, as in the genus $0$ case. Assume that the path $\gamma$ from $Q_{1}$ to $z$ that is involved in the iterated integral is a straight-line path. For $1\le i \le n-1$ one has 
$$\gamma^{*} w_{i} =(z-Q_{1}) \cdot (f(Q_{1}+t(z-Q_{1}) - P_{i_{1}}) - (f(Q_{1}+t(z-Q_{1}) - P_{i_{2}}) \textrm{ d}t =: F_{i}(t) \textrm{ d}t$$
For 
$$I(t_{k}) :=\int_{0}^{t_{k}} F_{j_{k+1}} \left ( \int_{0}^{t_{k-1}} F_{j_{k+2}}\cdots \right  ) \textrm{ d}t_{k-1}$$
one has

\begin{eqnarray*}
(\partial_{P_{j_{k_{1}}}} + \partial_{P_{j_{k_{2}}}}) \;  \textrm{Li}_{w_{j_{1}}\cdots w_{j_{k}} \cdots w_{j_{r}}}(Q_{1},Q_{2}) &=&(\partial_{P_{j_{k_{1}}}} + \partial_{P_{j_{k_{2}}}}) \;  \int_{t_{r}=0}^{1} F_{j_{1}}  \left ( \int_{t_{r-1}=0}^{t_{r}} F_{j_{2}} \left (\cdots \int_{t_{k}=0}^{t_{k+1}} F_{j_{k}}(t_{k}) \cdot I(t_{k}) \textrm{ d}t_{k} \right )  \textrm{ d}t_{r-1} \right ) \textrm{ d}t_{r}
\end{eqnarray*}
It follows that 
\begin{align*}
&  (z-Q_{1})^{-1}  \int_{t_{r}=0}^{1} F_{j_{1}} \left ( \int_{t_{r-1}=0}^{t_{r}} F_{j_{2}} \left ( \cdots \int_{0}^{t_{k+1}} - \partial_{t_{k}} F_{k}(t_{k}) \cdot I(t_{k}) \textrm{ d}t_{k} \right ) \textrm{ d}t_{r-1} \right )\textrm{ d}t_{r}\\
&\\
& =  -  \int_{t_{r}=0}^{1} F_{j_{1}}  \left (\int_{t_{r-1}=0}^{t_{r}} F_{j_{2}} \left (  \cdots \int_{0}^{t_{k+2}} \; \sum_{i=0}^{n-1} C_{j_{k+1},j_{k}}^{(i)} \cdot  F_{i}(t_{k+1}) \cdot I(t_{k+1}) \textrm{ d} t_{k+1}  \right ) \textrm{ d}t_{r-1} \right  ) \textrm{ d}t_{r} \\
& \\
 &  + \int_{t_{r}=0}^{1} F_{j_{1}} \left (\int_{t_{r-1}=0}^{t_{r}} F_{j_{2}} \cdots \int_{0}^{t_{k+1}}\sum_{i=0}^{n-1} \; C_{j_{k},j_{k-1}}^{(i)} \cdot  F_{i}(t_{k})  \cdot \left (\left ( \int_{0}^{t_{k}} F_{j_{k+1}}\cdots \right  ) \textrm{ d}t_{k-1} \right ) \textrm{ d}t_{k} \right )\textrm{ d}t_{r}\\
\end{align*}
Hence one obtains
\begin{eqnarray}
\label{differential-equation}
& (\partial_{P_{j_{k_{1}}}} + \partial_{P_{j_{k_{2}}}}) \; \textrm{Li}_{w_{j_{1}}\cdots w_{j_{k}} \cdots w_{j_{r}}}(z)  \nonumber \\
& =\\
&  \sum_{i=0}^{n-1}  \left ( C_{j_{k},j_{k+1}}^{(i)} \cdot  \textrm{Li}_{w_{j_{1}}\cdots \widehat{w}_{j_{k}} \widehat{w}_{j_{k+1}}  w_{i} \cdots w_{j_{r}}}(z)  -C_{j_{k-1},j_{k}}^{(i)} \cdot  \textrm{Li}_{w_{j_{1}}\cdots \widehat{w}_{j_{k-1}} \widehat{w}_{j_{k}} w_{i} \cdots w_{j_{r}}}(z) \right ) \nonumber
\end{eqnarray}
in complete analogy with the genus $0$ formula given in Equation (\ref{sphere-differential-equation}).
One can combine Equation (\ref{product-equation}) (for $(a,b)=(j_{k},j_{k-1})$ and for $(a,b)=(j_{k+1},j_{k})$) with Equation (\ref{differential-equation}) to calculate explicitly the action of $\partial_{P_{j_{k_{1}}}} +  \partial_{P_{j_{k_{2}}}}$ on the function $\textrm{Li}_{w_{j_{1}}\cdots w_{j_{k}} \cdots w_{j_{r}}}(z)$.

\hspace{0.2in}

\noindent
\textbf{Acknowledgements:} It is a great pleasure to thank Dmitry Fuchs and Andrey Levin for very helpful exchanges.

\Addresses


\begin{thebibliography}{99}
\bibitem[1]{BRO} F. Brown: Iterated integrals in quantum field theory, IHES (2009)
\bibitem[2]{BL} F. Brown, A. Levin: Multiple elliptic polylogarithms, arXiv:math/1110.6917
\bibitem[3]{BMMS} J. Broedel, C. Mafra, N.  Matthes, O. Schlotterer: Elliptic multiple zeta values and one-loop superstring amplitudes, JHEP \textbf{07} (2015)
\bibitem[4]{CHE} K. T. Chen: Iterated path integrals, Bull. Amer. Math. Soc. \textbf{83} (1977), 831 - 879 
\bibitem[5]{GON} A. Goncharov: Multiple polylogarithms and mixed Tate motives, arXiv:math/0103059
\bibitem[6]{HAI} R. Hain: Iterated integrals and algebraic cycles: Examples and prospects, Contemporary Trends in Algebraic Geometry and Algebraic Topology, Nankai Tracts in Mathematics \textbf{5}, World Scientific (2002)
\end{thebibliography}
\end{document}